\documentclass[12pt]{conm-p-l}

\usepackage{amssymb}
\usepackage{epsfig}
\usepackage{cite}

\newtheorem{theorem}{Theorem}[section]

\theoremstyle{definition}

\newtheorem{example}[theorem]{Example}

\theoremstyle{remark}

\numberwithin{equation}{section}

\newcommand{\cle}{\preceq}
\newcommand{\opl}{{\oplus}}

\newcommand{\rmin}{\mathbf{R}_{\min}}
\newcommand{\rmax}{\mathbf{R}_{\max}}
\newcommand{\suplim}{\sup\limits}
\newcommand{\sumlim}{\sum\limits}
\newcommand{\maxlim}{\max\limits}
\newcommand{\pd}[2]{\dfrac{\partial#1}{\partial#2}}
\newcommand{\CalB}{\mathcal{B}}

\newcommand{\cF}{{\mathcal F}}

\newcommand{\cN}{{\mathcal N}}

\newcommand{\0}{\mathbf{0}}
\newcommand{\1}{\mathbf{1}}
\newcommand{\cset}{\mathbf{C}}

\newcommand{\rset}{\mathbf{R}}
\newcommand{\maA}{\mathcal{A}}
\newcommand{\Log}{\mathop{\mathrm{Log}}}
\newcommand{\ovol}{\mathop{\mathrm{vol}}}

\def\C{\mathbf C}
\def\R{\mathbf R}

\def\cN{\mathcal N}

\begin{document}

\title{Tropical Mathematics, Classical Mechanics and Geometry}
\author{G.~L.~Litvinov}
\subjclass{Primary: 15A80, 46S19, 81Q20, 14M25, 16S80, 70H20, 14T05, 51P05, 52A20;
Secondary: 81S99, 52B70, 12K10, 46L65, 11K55, 28B10, 28A80, 28A25, 06F99, 16H99.}
\keywords{Tropical mathematics, idempotent
mathematics, classical mechanics, convex geometry, tropical
geometry, Newton polytopes.}

\thanks{This work is supported by the RFBR
grant 08--01--00601.}

\address{Independent University of Moscow,
Bol'shoi Vlasievskii per., 11, Moscow 119002, Russia}

\email{glitvinov@gmail.com}

\begin{abstract}
A very brief introduction to tropical and idempotent mathematics is presented. Applications to classical mechanics and geometry are especially examined.
\end{abstract}

\maketitle

\normalsize

\section{Introduction}

Tropical mathematics can be treated  as  a result of a
dequantization of the traditional mathematics as  the Planck
constant  tends to zero  taking imaginary values. This kind of
dequantization is known as the Maslov dequantization and it leads
to a mathematics over tropical algebras like the max-plus algebra.
The so-called idempotent dequantization is a generalization of the
Maslov dequantization. The idempotent dequantization leads to
mathematics over idempotent semirings (exact definitions see below
in sections 2 and 3). For  example,
the field of real or complex numbers can be treated as a quantum
object whereas idempotent semirings can be examined  as
"classical" or "semiclassical" objects  (a semiring is called
idempotent  if the  semiring addition is idempotent, i.e. $x
\oplus x = x$), see~\cite{Li2005,LiMa95,LiMa96,LiMa98}.

Tropical algebras are idempotent semirings (and semifields). Thus tropical mathematics is a part of idempotent mathematics. Tropical
algebraic geometry can be treated as a result of the Maslov
dequantization applied to the traditional algebraic geometry (O.
Viro, G. Mikhalkin), see,
e.g.,~\cite{Iten2007,Mi2005,Mi2006,Vir2000,Vir2002,Vir2008}. There
are interesting relations and applications to the traditional
convex geometry.

    In the spirit of N.~Bohr's correspondence principle there
is a (heuristic)  correspondence  between important, useful, and
interesting  constructions  and  results over fields and similar
results  over  idempotent  semirings.  A systematic application of
this correspondence  principle  leads  to  a variety  of theoretical
and applied results~\cite{Li2005,LiMa95,LiMa96,LiMa98,LiMa2005}, see
Fig.~1.

\begin{figure}
\centering
\epsfig{file=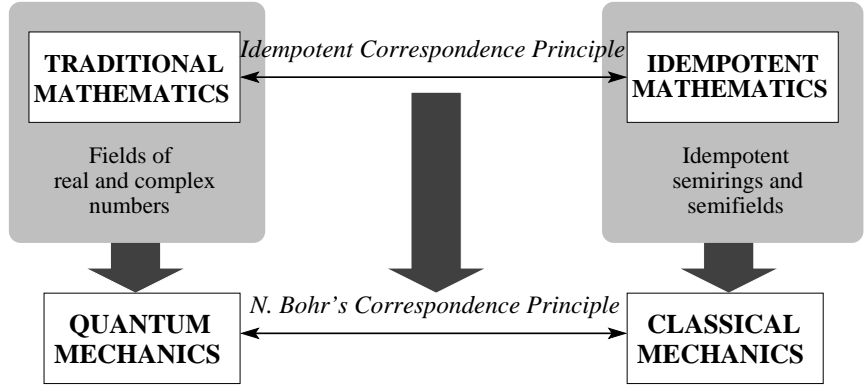,width=0.9\linewidth} \caption{Relations
between idempotent and traditional mathematics.}
\end{figure}

The history of the subject is discussed, e.g., in~\cite{Li2005}.
There is a large list of references.

\section{The Maslov dequantization}

Let $\R$ and $\C$ be the fields of real and complex numbers. The
so-called max-plus algebra $\R_{\max}= \R\cup\{-\infty\}$ is
defined by the operations $x\oplus y=\max\{x, y\}$ and $x\odot y=
x+y$.

The max-plus algebra can be treated as a result of the {\it Maslov
dequantization} of the semifield $\R_+$ of all nonnegative
numbers with the usual arithmetics. The change of variables
\begin{eqnarray*}
x\mapsto u=h\log x,
\end{eqnarray*}
where $h>0$, defines a map $\Phi_h\colon \R_+\to
\R\cup\{-\infty\}$, see Fig.~2. Let the addition and multiplication
operations be mapped from \markboth{G.L. Litvinov}{Tropical Mathematics, Classical Mechanics, and Geometry} $\R_+$ to $\R\cup\{-\infty\}$ by $\Phi_h$, i.e.\
let
\begin{eqnarray*}
u\oplus_h v = h \log({\mbox{exp}}(u/h)+{\mbox{exp}}(v/h)),\quad u\odot v= u+ v,\\
\mathbf{0}=-\infty = \Phi_h(0),\quad \mathbf{1}= 0 = \Phi_h(1).
\end{eqnarray*}

\begin{figure}
\noindent\epsfig{file=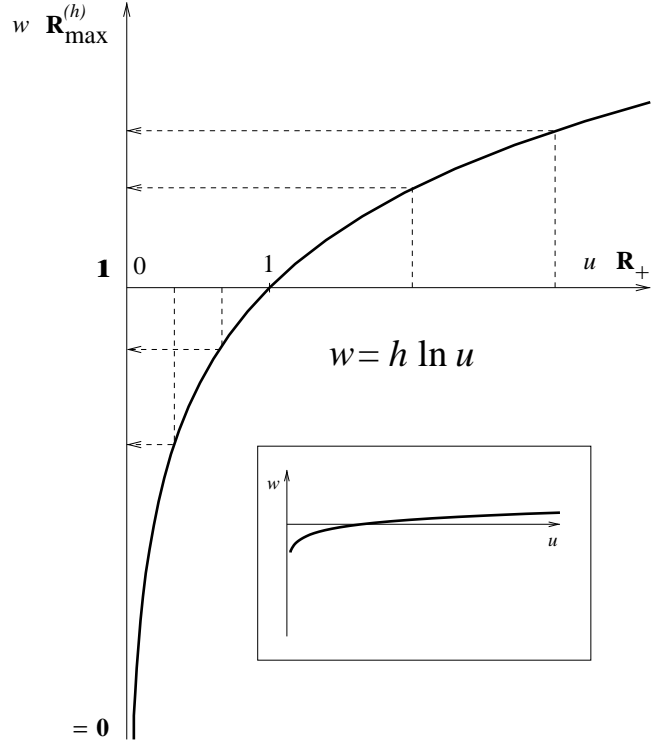,width=0.9\linewidth}
\vskip -3cm
\caption{Deformation of $\R_+$ to $\R^{(h)}$. Inset: the same for a
small value of $h$.}
\end{figure}

It can easily be checked that $u\oplus_h v\to \max\{u, v\}$ as
$h\to 0$. Thus we get the semifield $\R_{\max}$ (i.e.\ the
max-plus algebra) with zero $\mathbf{0}= -\infty$ and unit
$\mathbf{1}=0$ as a result of this deformation of the algebraic
structure~in~$\R_+$.

The semifield $\R_{\max}$ is a typical example of an {\it
 idempotent semiring}; this is a semiring with idempotent addition, i.e.,
 $x\oplus x = x$ for arbitrary element
 $x$ of this semiring.

 The semifield $\R_{\max}$ is also called a \emph{tropical
 algebra}. The semifield $\R^{(h)}=\Phi_h(\R_+)$ with operations
 $\oplus_h$ and $\odot$ (i.e.$+$) is called a \emph{subtropical
 algebra}.

 The semifield $\R_{\min}=\R\cup\{+\infty\}$
with operations $\oplus={\min}$ and $\odot=+$
$(\mathbf{0}=+\infty, \mathbf{1}=0)$ is isomorphic to $\R_{\max}$.

The analogy with quantization is obvious; the parameter $h$ plays
the role of the Planck constant. The map $x\mapsto|x|$ and the
Maslov dequantization for $\R_+$ give us a natural transition from
the field $\C$ (or $\R$) to the max-plus algebra $\R_{\max}$. {\it
We will also call this transition the Maslov dequantization}. In fact
the Maslov dequantization corresponds to the usual Schr\"odinger
dequantization but for  imaginary values of the Planck constant (see below).
The transition from numerical fields to the max-plus algebra
$\R_{\max}$ (or similar semifields) in mathematical constructions
and results generates the so called {\it tropical mathematics}.
The so-called {\it idempotent dequantization} is a generalization
of the Maslov dequantization; this is the transition from basic fields to idempotent semirings in mathematical constructions and results without any deformation. The idempotent dequantization generates
the so-called \emph{idempotent mathematics}, i.e. mathematics over
idempotent semifields and semirings.

{\bf Remark.} The term 'tropical' appeared in~\cite{Sim88} for a
discrete version of the max-plus algebra (as a suggestion of
Christian Choffrut). On the other hand V.~P.~Maslov used this term in
80s in his talks and works on economical applications of his
idempotent analysis (related to colonial politics). For the most
part of modern authors, 'tropical' means 'over $\R_{\max}$ (or
$\R_{\min}$)' and tropical algebras are $\R_{\max}$ and $\R_{\min}$.
The terms 'max-plus', 'max-algebra' and 'min-plus' are often used in
the same sense.
\medskip

\section{Semirings and semifields}

Consider a set $S$ equipped with two algebraic operations: {\it
addition} $\oplus$ and {\it multiplication} $\odot$. It is a {\it
semiring} if the following conditions are satisfied:
\begin{itemize}
\item the addition $\oplus$ and the multiplication $\odot$ are
associative; \item the addition $\oplus$ is commutative; \item the
multiplication $\odot$ is distributive with respect to the
addition $\oplus$:
\[x\odot(y\oplus z)=(x\odot y)\oplus(x\odot z)\]
and
\[(x\oplus y)\odot z=(x\odot z)\oplus(y\odot z)\]
for all $x,y,z\in S$.
\end{itemize}
A {\it unity} (we suppose that it exists) of a semiring $S$ is an element $\1\in S$ such that
$\1\odot x=x\odot\1=x$ for all $x\in S$. A {\it zero} (if it exists) of a
semiring $S$ is an element $\0\in S$ such that $\0\neq\1$ and
$\0\oplus x=x$, $\0\odot x=x\odot \0=\0$ for all $x\in S$. A
semiring $S$ is called an {\it idempotent semiring} if $x\oplus
x=x$ for all $x\in S$. A semiring $S$ with neutral element
$\1$ is called a {\it semifield} if every nonzero element of
$S$ is invertible with respect to the multiplication. The theory
of semirings and semifields is treated, e.g., in~\cite{Gol99}.
\medskip

\section{Idempotent analysis}

Idempotent analysis deals with functions taking their values in
an idempotent semiring and the corresponding function spaces.
Idempotent analysis was initially constructed by V.~P.~Maslov and
his collaborators and then developed by many authors. The subject
and applications are presented in the book of V.~N.~Kolokoltsov and
V.~P.~Maslov~\cite{KoMa97} (a version of this book in Russian was
published in 1994).

Let $S$ be an arbitrary semiring with idempotent addition $\oplus$
(which is always assumed to be commutative), multiplication
$\odot$, and unit $\1$. The set $S$ is supplied with
the {\it standard partial order\/}~$\cle$: by definition, $a \cle
b$ if and only if $a \oplus b = b$. If $S$ contains a zero element
$\0$, then all elements of $S$ are
nonnegative: $\0 \cle$ $a$ for all $a \in S$. Due to the existence
of this order, idempotent analysis is closely related to the
lattice theory, theory of vector lattices, and theory of ordered
spaces. Moreover, this partial order allows to model a number of
basic ``topological'' concepts and results of idempotent analysis
at the purely algebraic level; this line of reasoning was examined
systematically in~\cite{Li2005}--~\cite{LiSh2007}
and~\cite{CoGaQu2004}.

Calculus deals mainly with functions whose values are numbers. The
idempotent analog of a numerical function is a map $X \to S$,
where $X$ is an arbitrary set and $S$ is an idempotent semiring.
Functions with values in $S$ can be added, multiplied by each
other, and multiplied by elements of $S$ pointwise.

The idempotent analog of a linear functional space is a set of
$S$-valued functions that is closed under addition of functions
and multiplication of functions by elements of $S$, or an
$S$-semimodule. Consider, e.g., the $S$-semimodule $B(X, S)$ of
all functions $X \to S$ that are bounded in the sense of the
standard order on $S$.

If $S = \rmax$, then the idempotent analog of integration is
defined by the formula
$$
I(\varphi) = \int_X^{\oplus} \varphi (x)\, dx     = \sup_{x\in X}
\varphi (x),\eqno{(1)}
$$
where $\varphi \in B(X, S)$. Indeed, a Riemann sum of the form
$\sumlim_i \varphi(x_i) \cdot \sigma_i$ corresponds to the
expression $\bigoplus\limits_i \varphi(x_i) \odot \sigma_i =
\maxlim_i \{\varphi(x_i) + \sigma_i\}$, which tends to the
right-hand side of~(1) as $\sigma_i \to 0$. Of course, this is a
purely heuristic argument.

Formula~(1) defines the \emph{idempotent} (or \emph{Maslov})
\emph{integral} not only for functions taking values in $\rmax$,
but also in the general case when any of bounded (from above)
subsets of~$S$ has the least upper bound.

An \emph{idempotent} (or \emph{Maslov}) \emph{measure} on $X$ is
defined by the formula $m_{\psi}(Y) = \suplim_{x \in Y} \psi(x)$, where $\psi
\in B(X,S)$ is a fixed function. The integral with respect to this measure is defined
by the formula
$$
   I_{\psi}(\varphi)
    = \int^{\oplus}_X \varphi(x)\, dm_{\psi}
    = \int_X^{\oplus} \varphi(x) \odot \psi(x)\, dx
    = \sup_{x\in X} (\varphi (x) \odot \psi(x)).
    \eqno{(2)}
$$

Obviously, if $S = \rmin$, then the standard order is
opposite to the conventional order $\le$, so in this case
equation~(2) assumes the form
$$
   \int^{\oplus}_X \varphi(x)\, dm_{\psi}
    = \int_X^{\oplus} \varphi(x) \odot \psi(x)\, dx
    = \inf_{x\in X} (\varphi (x) \odot \psi(x)),
$$
where $\inf$ is understood in the sense of the conventional order
$\le$.
\medskip

\section{The superposition principle and linear problems}

Basic equations of quantum theory are linear; this is the
superposition principle in quantum mechanics. The Hamilton--Jacobi
equation, the basic equation of classical mechanics, is nonlinear
in the conventional sense. However, it is linear over the
semirings $\rmax$ and $\rmin$. Similarly, different versions of
the Bellman equation, the basic equation of optimization theory,
are linear over suitable idempotent semirings. This is
V.~P.~Maslov's idempotent superposition principle, see
\cite{Mas86,Mas87a,Mas87b}. For instance, the finite-dimensional
stationary Bellman equation can be written in the form $X = H
\odot X \oplus F$, where $X$, $H$, $F$ are matrices with
coefficients in an idempotent semiring $S$ and the unknown matrix
$X$ is determined by $H$ and $F$~\cite{Ca71, Ca79, BaCoOlQu92, Cu79, Cu95, GoMi79, GoMi2001}. In
particular, standard problems of dynamic programming and the
well-known shortest path problem correspond to the cases $S =
\rmax$ and $S =\rmin$, respectively. It is known that principal
optimization algorithms for finite graphs correspond to standard
methods for solving systems of linear equations of this type
(i.e., over semirings). Specifically, Bellman's shortest path
algorithm corresponds to a version of Jacobi's algorithm, Ford's
algorithm corresponds to the Gauss--Seidel iterative scheme,
etc.~\cite{Ca71, Ca79}.

The linearity of the Hamilton--Jacobi equation over $\rmin$ and
$\rmax$, which is the result of the Maslov dequantization of the
Schr{\"o}\-din\-ger equation, is closely related to the
(conventional) linearity of the Schr{\"o}\-din\-ger equation and
can be deduced from this linearity. Thus, it is possible to borrow
standard ideas and methods of linear analysis and apply them to a
new area.

Consider a classical dynamical system specified by the Hamiltonian
$$
   H = H(p,x) = \sum_{i=1}^N \frac{p^2_i}{2m_i} + V(x),
$$
where $x = (x_1, \dots, x_N)$ are generalized coordinates, $p =
(p_1, \dots, p_N)$ are generalized momenta, $m_i$ are generalized
masses, and $V(x)$ is the potential. In this case the Lagrangian
$L(x, \dot x, t)$ has the form
$$
   L(x, \dot x, t)
    = \sum^N_{i=1} m_i \frac{\dot x_i^2}2 - V(x),
$$
where $\dot x = (\dot x_1, \dots, \dot x_N)$, $\dot x_i = dx_i /
dt$. The value function $S(x,t)$ of the action functional has the
form
$$
   S = \int^t_{t_0} L(x(t), \dot x(t), t)\, dt,
$$
where the integration is performed along the factual trajectory of
the system.  The classical equations of motion are derived as the
stationarity conditions for the action functional (the Hamilton
principle, or the least action principle).

For fixed values of $t$ and $t_0$ and arbitrary trajectories
$x(t)$, the action functional $S=S(x(t))$ can be considered as a
function taking the set of curves (trajectories) to the set of
real numbers which can be treated as elements of  $\rmin$. In this
case the minimum of the action functional can be viewed as the
Maslov integral of this function over the set of trajectories or
an idempotent analog of the Euclidean version of the Feynman path
integral. The minimum of the action functional corresponds to the
maximum of $e^{-S}$, i.e. idempotent integral
$\int^{\oplus}_{\{paths\}} e^{-S(x(t))} D\{x(t)\}$ with respect to
the max-plus algebra $\rset_{\max}$. Thus the least action
principle can be considered as an idempotent version of the
well-known Feynman approach to quantum mechanics.  The
representation of a solution to the Schr{\"o}\-din\-ger equation
in terms of the Feynman integral corresponds to the
Lax--Ole\u{\i}nik solution formula for the Hamilton--Jacobi
equation.

Since $\partial S/\partial x_i = p_i$, $\partial S/\partial t =
-H(p,x)$, the following Hamilton--Jacobi equation holds:
$$
   \pd{S}{t} + H \left(\pd{S}{x_i}, x_i\right)= 0.\eqno{(3)}
$$

Quantization leads to the Schr\"odinger equation
$$
   -\frac{\hbar}i \pd{\psi}{t}= \widehat H \psi = H(\hat p_i, \hat x_i)\psi,
   \eqno{(4)}
$$
where $\psi = \psi(x,t)$ is the wave function, i.e., a
time-dependent element of the Hilbert space $L^2(\rset^N)$, and
$\widehat H$ is the energy operator obtained by substitution of
the momentum operators $\widehat p_i = {\hbar \over i}{\partial
\over \partial x_i}$ and the coordinate operators $\widehat x_i
\colon \psi \mapsto x_i\psi$ for the variables $p_i$ and $x_i$ in
the Hamiltonian function, respectively. This equation is linear in
the conventional sense (the quantum superposition principle). The
standard procedure of limit transition from the Schr\"odinger
equation to the Hamilton--Jacobi equation is to use the following
ansatz for the wave function:  $\psi(x,t) = a(x,t)
e^{iS(x,t)/\hbar}$, and to keep only the leading order as $\hbar
\to 0$ (the `semiclassical' limit).

Instead of doing this, we switch to imaginary values of the Planck
constant $\hbar$ by the substitution $h = i\hbar$, assuming $h >
0$. Thus the Schr\"odinger equation~(4) turns to an analog of the
heat equation:
$$
   h\pd{u}{t} = H\left(-h\frac{\partial}{\partial x_i}, \hat x_i\right) u,
   \eqno{(5)}
$$
where the real-valued function $u$ corresponds to the wave
function $\psi$. A similar idea (the switch to imaginary time) is
used in the Euclidean quantum field theory; let us remember that
time and energy are dual quantities.

Linearity of equation~(4) implies linearity of equation~(5). Thus
if $u_1$ and $u_2$ are solutions of~(5), then so is their linear
combination
$$
   u = \lambda_1 u_1 + \lambda_2 u_2.\eqno{(6)}
$$

Let $S = h \ln u$ or $u = e^{S/h}$ as in Section 2 above. It can
easily be checked that equation~(5) thus turns to
$$
   \pd{S}{t}= V(x) + \sum^N_{i=1} \frac1{2m_i}\left(\pd{S}{x_i}\right)^2
    + h\sum^n_{i=1}\frac1{2m_i}\frac{\partial^2 S}{\partial x^2_i}.
   \eqno{(7)}
$$
Thus we have a transition from (4) to (7) by means of the change of
variables $\psi = e^{S/h}$. Note that $|\psi| = e^{ReS/h}$ , where
Re$S$ is the real part  of $S$. Now let us consider $S$ as a real
variable. The equation (7) is nonlinear in the conventional sense.
However, if $S_1$ and $S_2$ are its solutions, then so is the
function
$$
   S = \lambda_1 \odot S_1 \opl_h \lambda_2\odot S_2
$$
obtained from~(6) by means of our substitution $S = h \ln u$. Here
the generalized multiplication $\odot$ coincides with the ordinary
addition and the generalized addition $\opl_h$ is the image of the
conventional addition under the above change of variables.  As $h
\to 0$, we obtain the operations of the idempotent semiring
$\rmax$, i.e., $\oplus = \max$ and $\odot = +$, and equation~(7)
turns to the Hamilton--Jacobi equation~(3), since the third term
in the right-hand side of equation~(7) vanishes.

Thus it is natural to consider the limit function $S = \lambda_1
\odot S_1 \oplus \lambda_2 \odot S_2$ as a solution of the
Hamilton--Jacobi equation and to expect that this equation can be
treated as linear over $\rmax$. This argument (clearly, a
heuristic one) can be extended to equations of a more general
form. For a rigorous treatment of (semiring) linearity for these
equations see, e.g., \cite{KoMa97,LiMa2005,Roub}. Notice that if
$h$ is changed to $-h$, then we have that the resulting
Hamilton--Jacobi equation is linear over $\rmin$.

The idempotent superposition principle indicates that there exist
important nonlinear (in the traditional sense) problems that are
linear over idempotent semirings. The idempotent linear functional
analysis (see below) is a natural tool for investigation of those
nonlinear infinite-dimensional problems that possess this
property.

\section{Convolution and the Fourier--Legendre transform}

Let $G$ be a group. Then the space $\CalB(G, \rset_{\max})$ of all
bounded functions $G\to\rset_{\max}$ (see above) is an idempotent
semiring with respect to the following analog $\circledast$ of the
usual convolution:
$$
   (\varphi(x)\circledast\psi)(g)=
    = \int_G^{\oplus} \varphi (x)\odot\psi(x^{-1}\cdot g)\, dx=
\sup_{x\in G}(\varphi(x)+\psi(x^{-1}\cdot g)).
$$
Of course, it is possible to consider other ``function spaces''
(and other basic semirings instead of $\rset_{\max}$).

Let $G=\rset^n$, where $\rset^n$ is considered as a topological
group with respect to the vector addition. The conventional
Fourier--Laplace transform is defined as
$$
   \varphi(x) \mapsto \tilde{\varphi}(\xi)
    = \int_G e^{i\xi \cdot x} \varphi (x)\, dx,\eqno{(8)}
$$
where $e^{i\xi \cdot x}$ is a character of the group $G$, i.e., a
solution of the following functional equation:
$$
   f(x + y) = f(x)f(y).
$$
The idempotent analog of this equation is
$$
   f(x + y) = f(x) \odot f(y) = f(x) + f(y),
$$
so ``continuous idempotent characters'' are linear functionals of
the form $x \mapsto \xi \cdot x = \xi_1 x_1 + \dots + \xi_n x_n$.
As a result, the transform in~(8) assumes the form
$$
   \varphi(x) \mapsto \tilde{\varphi}(\xi)
    = \int_G^\oplus \xi \cdot x \odot \varphi (x)\, dx
   = \sup_{x\in G} (\xi \cdot x + \varphi (x)).\eqno{(9)}
$$
The transform in~(9) is nothing but the {\it Legendre transform\/}
(up to some notation) \cite{Mas87b}; transforms of this kind
establish the correspondence between the Lagrangian and the
Hamiltonian formulations of classical mechanics. The Legendre
transform generates an idempotent version of harmonic analysis for
the space of convex functions, see, e.g., \cite{MaTi2003}.

Of course, this construction can be generalized to different
classes of groups and semirings. Transformations of this type
convert the generalized convolution $\circledast$ to the pointwise
(generalized) multiplication and possess analogs of some important
properties of the usual Fourier transform.

The examples discussed in this sections can be treated as
fragments of an idempotent version of the representation theory,
see, e.g., \cite{LiMaSh2002}. In particular, ``idempotent''
representations of groups can be examined as representations of
the corresponding convolution semirings (i.e. idempotent group
semirings) in semimodules.
\medskip

\section{Idempotent functional analysis}

Many other idempotent analogs may be given, in particular, for
basic constructions and theorems of functional analysis.
Idempotent functional analysis is an abstract version of
idempotent analysis. For the sake of simplicity take $S=\rmax$ and
let $X$ be an arbitrary set. The idempotent integration can be
defined by the formula (1), see above. The functional $I(\varphi)$
is linear over $S$ and its values correspond to limiting values of
the corresponding analogs of Lebesgue (or Riemann) sums. An
idempotent scalar product of functions $\varphi$ and $\psi$ is
defined by the formula
$$
\langle\varphi,\psi\rangle = \int^{\oplus}_X
\varphi(x)\odot\psi(x)\, dx = \sup_{x\in
X}(\varphi(x)\odot\psi(x)).
$$
So it is natural to construct idempotent analogs of integral
operators in the form
$$
\varphi(y) \mapsto (K\varphi)(x) = \int^{\oplus}_Y K(x,y)\odot
\varphi(y)\, dy = \sup_{y\in Y}\{K(x,y)+\varphi(y)\},\eqno(10)
$$
where $\varphi(y)$ is an element of a space of functions defined
on a set $Y$, and $K(x,y)$ is an $S$-valued function on $X\times
Y$. Of course, expressions of this type are standard in
optimization problems.\medskip

Recall that the definitions and constructions described above can
be extended to the case of idempotent semirings which are
conditionally complete in the sense of the standard order. Using
the Maslov integration, one can construct various function spaces
as well as idempotent versions of the theory of generalized
functions (distributions). For some concrete idempotent function
spaces it was proved that every `good' linear operator (in the
idempotent sense) can be presented in the form (10); this is an
idempotent version of the kernel theorem of L.~Schwartz; results
of this type were proved by V.~N.~Kolokoltsov, P.~S.~Dudnikov and
S.~N.~Samborski\u\i, I.~Singer, M.~A.~Shubin and others. So every
`good' linear functional can be presented in the form
$\varphi\mapsto\langle\varphi,\psi\rangle$, where
$\langle,\rangle$ is an idempotent scalar product.\medskip

In the framework of idempotent functional analysis results of this
type can be proved in a very general situation. In \cite{LiMaSh98,
LiMaSh99,LiMaSh2001,LiMaSh2002,LiSh2002,LiSh2007} an algebraic
version of the idempotent functional analysis is developed; this
means that basic (topological) notions and results are simulated
in purely algebraic terms (see below). The treatment covers the subject  from
basic concepts and results (e.g., idempotent analogs of the
well-known theorems of Hahn-Banach, Riesz, and Riesz-Fisher) to
idempotent analogs of A.~Grothendieck's concepts and results on
topological tensor products, nuclear spaces and operators.
Abstract idempotent versions of the kernel theorem are formulated. Note that
the transition from the usual theory to idempotent functional
analysis may be very nontrivial; for example, there are many
non-isomorphic idempotent Hilbert spaces. Important results on
idempotent functional analysis (duality and separation theorems)
were obtained by G.~Cohen, S.~Gaubert, and J.-P.~Quadrat.
Idempotent functional analysis has received much attention in the
last years, see,
e.g.,~\cite{AkGaKo2005, CoGaQu2004, GoMi79, GoMi2001, Gun98a, MaSa92, Shub92},~\cite{KoMa97}--~\cite{LiSh2007} and works
cited in~\cite{Li2005}. Elements of "tropical" functional analysis are presented in~\cite{KoMa97}. 

\medskip

\section{The dequantization transform, convex geometry and the Newton polytopes}

Let $X$ be a topological space. For functions $f(x)$ defined on
$X$ we shall say that a certain property is valid {\it almost
everywhere} (a.e.) if it is valid for all elements $x$ of an open
dense subset of $X$. Suppose $X$ is $\C^n$ or $\R^n$; denote by
$\R^n_+$ the set $x=\{\,(x_1, \dots, x_n)\in X \mid x_i\geq 0$ for
$i = 1, 2, \dots, n$.
 For $x= (x_1, \dots, x_n) \in X$ we set
${\mbox{exp}}(x) = ({\mbox{exp}}(x_1), \dots, {\mbox{exp}}(x_n))$;
so if $x\in\R^n$, then ${\mbox{exp}}(x)\in \R^n_+$.

Denote by $\cF(\C^n)$ the set of all functions defined and
continuous on an open dense subset $U\subset \C^n$ such that
$U\supset \R^n_+$. It is clear that $\cF(\C^n)$ is a ring (and an
algebra over $\C$) with respect to the usual addition and
multiplications of functions.

For $f\in \cF(\C^n)$ let us define the function $\hat f_h$ by the
following formula:
$$
\hat f_h(x) = h \log|f({\mbox{exp}}(x/h))|,
\eqno(11)
$$
where $h$ is a (small) real positive parameter and $x\in\R^n$. Set
$$
\label{e:hatfx} \hat f(x) = \lim_{h\to +0} \hat f_h (x),
\eqno(12)
$$
if the right-hand side of (12) exists almost everywhere.

We shall say that the function $\hat f(x)$ is a {\it
dequantization} of the function $f(x)$ and the map $f(x)\mapsto
\hat f(x)$ is a {\it dequantization transform}. By construction,
$\hat f_h(x)$ and $\hat f(x)$ can be treated as functions taking
their values in $\R_{\max}$. Note that in fact $\hat f_h(x)$ and
$\hat f(x)$
 depend on the restriction of $f$ to $\R_+^n $
only; so in fact the dequantization transform is constructed for
functions defined on $\R^n_+$ only. It is clear that the
dequantization transform is generated by the Maslov dequantization
and the map $x\mapsto |x|$.

Of course, similar definitions can be given for functions defined
on $\R^n$ and $\R_+^n$. If $s=1/h$, then we have the following
version of (11) and (12):
$$
\hat f(x) = \lim_{s\to \infty} (1/s) \log|f(e^{sx})|.
\eqno(12')
$$

Denote by $\partial \hat f$ the subdifferential of the function
$\hat f$ at the origin.

If $f$ is a polynomial  we have
$$
\partial \hat f = \{\, v\in \R^n\mid (v, x) \le \hat f(x)\
\forall x\in \R^n\,\}.
$$
It is well known that all the convex compact subsets in $\R^n$
form an idempotent semiring $\mathcal{S}$ with respect to the
Minkowski operations: for $\alpha, \beta \in \mathcal{S}$ the sum
$\alpha\oplus \beta$ is the convex hull of the union $\alpha\cup
\beta$; the product $\alpha\odot \beta$ is defined in the
following way: $\alpha\odot \beta = \{\, x\mid x = a+b$, where
$a\in \alpha, b\in \beta$, see Fig.~3. In fact $\mathcal{S}$ is an
idempotent linear space over $\R_{\max}$.

Of course, the Newton polytopes of polynomials in $n$ variables
form a subsemiring
$\mathcal{N}$ in $\mathcal{S}$. If $f$, $g$ are polynomials, then
$\partial(\widehat{fg}) = \partial\hat f\odot\partial\widehat g$;
moreover, if $f$ and $g$ are ``in general position'', then
$\partial(\widehat{f+g}) = \partial\hat f\oplus\partial\widehat
g$. For the semiring of all polynomials with nonnegative
coefficients the dequantization transform is a homomorphism of
this ``traditional'' semiring to the idempotent semiring
$\mathcal{N}$.

\begin{figure}
\centering \epsfig{file=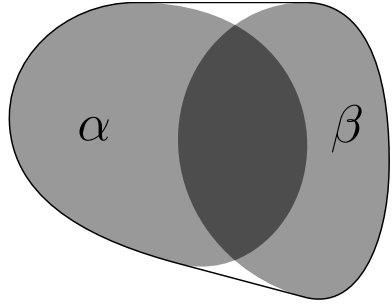,width=0.4\linewidth}
\caption{Algebra of convex subsets.}
\end{figure}

\begin{theorem}
If $f$ is a polynomial, then the subdifferential $\partial\hat f$
of $\hat f$ at the origin coincides with the Newton polytope of
$f$.  For the semiring of polynomials with nonnegative
coefficients, the transform $f\mapsto\partial\hat f$ is a
homomorphism of this semiring to the semiring of convex polytopes
with respect to the Minkowski operations (see above).
\end{theorem}

Using the dequantization transform it is possible to generalize
this result to a wide class of functions and convex
sets, see~\cite{LiSh2005}.

\medskip

\section{Dequantization of set functions and measures on metric spaces}
\medskip

The following results are presented in~\cite{LiSh}.

\begin{example} Let $M$ be a metric space, $S$ its arbitrary
subset with a compact closure. It is well-known that a Euclidean
$d$-dimensional ball $B_{\rho}$ of radius $\rho$ has volume
$$
\ovol\nolimits_d(B_{\rho})=\frac{\Gamma(1/2)^d}{\Gamma(1+d/2)}\rho^d,
$$
where $d$ is a natural parameter. By means of this formula it is
possible to define a volume of $B_{\rho}$ for any {\it real} $d$.
Cover $S$ by a finite number of balls of radii $\rho_m$. Set
$$
v_d(S):=\lim_{\rho\to 0} \inf_{\rho_m<\rho} \sum_m
\ovol\nolimits_d(B_{\rho_m}).
$$
Then there exists a number $D$ such that $v_d(S)=0$ for $d>D$ and
$v_d(S)=\infty$ for $d<D$. This number $D$ is called the {\it
Hausdorff-Besicovich dimension} (or {\it HB-dimension}) of $S$,
see, e.g.,~\cite{Ma2005}. Note that a set of non-integral
HB-dimension is called a fractal in the sense of B.~Mandelbrot.
\end{example}

\begin{theorem}
Denote by $\cN_{\rho}(S)$ the minimal number of balls of radius
$\rho$ covering $S$. Then
$$
D(S)=\mathop{\underline{\lim}}\limits_{\rho\to +0} \log_{\rho}
(\cN_{\rho}(S)^{-1}),
$$
where $D(S)$ is the HB-dimension of $S$. Set $\rho=e^{-s}$, then
$$
D(S)=\mathop{\underline{\lim}}\limits_{s\to +\infty} (1/s) \cdot
\log \cN_{exp(-s)}(S).
$$
So the HB-dimension $D(S)$ can be treated as a result of a
dequantization of the set function $\cN_{\rho}(S)$.
\end{theorem}

\begin{example} Let $\mu$ be a set function on $M$ (e.g., a
probability measure) and suppose that $\mu(B_{\rho})<\infty$ for
every ball $B_{\rho}$. Let $B_{x,\rho}$ be a ball of radius $\rho$
having the point $x\in M$ as its center. Then define
$\mu_x(\rho):=\mu(B_{x,\rho})$ and let $\rho=e^{-s}$ and
$$
D_{x,\mu}:=\mathop{\underline{\lim}}\limits_{s\to +\infty}
-(1/s)\cdot\log (|\mu_x(e^{-s})|).
$$
This number could be treated as a dimension of $M$ at the point
$x$ with respect to the set function $\mu$. So this dimension is a
result of a dequantization of the function $\mu_x(\rho)$, where
$x$ is fixed. There are many dequantization procedures of this
type in different mathematical areas. In particular, V.P.~Maslov's
negative dimension (see~\cite{Mas2007}) can be treated similarly.
\end{example}
\medskip

\section{Dequantization of geometry}

An idempotent version of real algebraic geometry was discovered in
the report of O.~Viro for the Barcelona Congress~\cite{Vir2000}.
Starting from the idempotent correspondence principle O.~Viro
constructed a piecewise-linear geometry of polyhedra of a special
kind in finite dimensional Euclidean spaces as a result of the
Maslov dequantization of real algebraic geometry. He indicated
important applications in real algebraic geometry (e.g.,
 in the framework of Hilbert's 16th problem for
constructing real algebraic varieties with prescribed properties
and parameters) and relations to complex algebraic geometry and
amoebas in the sense of I.~M.~Gelfand, M.~M.~Kapranov, and
A.~V.~Zelevinsky, see~\cite{GeKaZe94,Vir2002}. Then complex
algebraic geometry was dequantized by G.~Mikhalkin and the result
turned out to be the same; this new `idempotent' (or asymptotic)
geometry is now often called the {\it tropical algebraic
geometry}, see,
e.g.,~\cite{Iten2007,LiMa2005,LiSer2007,LiSer2009,Mi2005,Mi2006}.

There is a natural relation between the Maslov dequantization and
amoebas.

Suppose $({\cset}^*)^n$ is a complex torus, where ${\cset}^* =
{\cset}\backslash \{0\}$ is the group of nonzero complex numbers
under multiplication.  For
 $z = (z_1, \dots, z_n)\in
(\cset^*)^n$ and a positive real number $h$ denote by $\Log_h(z) =
h\log(|z|)$ the element
\[(h\log |z_1|, h\log |z_2|, \dots,
h\log|z_n|) \in \rset^n.\] Suppose $V\subset (\cset^*)^n$ is a
complex algebraic variety; denote by $\maA_h(V)$ the set
$\Log_h(V)$. If $h=1$, then the set $\maA(V) = \maA_1(V)$ is
called the {\it amoeba} of $V$; the amoeba $\maA(V)$ is a closed
subset of $\rset^n$ with a non-empty complement. Note that this
construction depends on our coordinate system.

For the sake of simplicity suppose $V$ is a hypersurface
in~$(\cset^*)^n$ defined by a polynomial~$f$; then there is a
deformation $h\mapsto f_h$ of this polynomial generated by the
Maslov dequantization and $f_h = f$ for $h = 1$. Let $V_h\subset
({\cset}^*)^n$ be the zero set of $f_h$ and set $\maA_h (V_h) =
{\Log}_h (V_h)$. Then
 there exists a tropical variety
$\mathit{Tro}(V)$ such that the subsets $\maA_h(V_h)\subset
\rset^n$ tend to $\mathit{Tro}(V)$ in the Hausdorff metric as
$h\to 0$. The tropical variety $\mathit{Tro}(V)$ is a result of a
deformation of the amoeba $\maA(V)$ and the Maslov dequantization
of the variety $V$. The set $\mathit{Tro}(V)$ is called the {\it
skeleton} of $\maA(V)$.

\begin{figure}
\noindent\epsfig{file=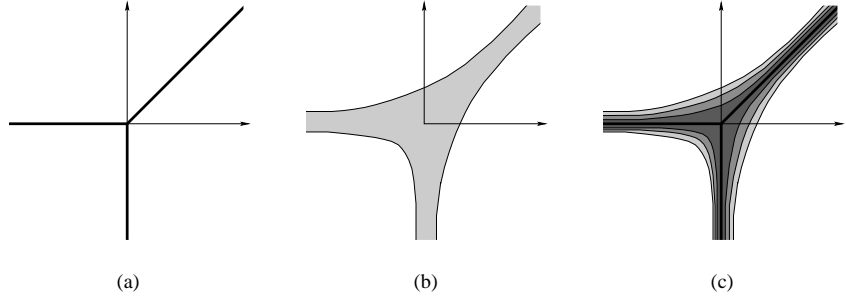,width=0.9\linewidth}
\caption{Tropical line and deformations of an amoeba.}
\end{figure}

\begin{example}  For the line $V = \{\, (x, y)\in ({\cset}^*)^2
\mid x + y + 1 = 0\,\}$ the piecewise-linear graph
$\mathit{Tro}(V)$ is a tropical line, see Fig.4(a). The amoeba
$\maA(V)$ is represented in Fig.4(b), while Fig.4(c) demonstrates
the corresponding deformation of the amoeba.
\end{example}

\medskip

{\bf Acknowledgments}.  The author is sincerely grateful to V.~N.~Kolokoltsov,
V.~P.~Maslov, G.~B.~Shpiz, S.~N.~Sergeev, and A.~N.~Sobolevski{{\u\i}} for valuable suggestions, help and support.
\medskip

\end{document}